\documentclass[12pt]{article}
\usepackage{amsfonts}
\newtheorem{theorem}{Theorem}
\newtheorem{prop}[theorem]{Proposition}
\newtheorem{unnumber}{}

\newtheorem{prepf}[unnumber]{Proof}
\newenvironment{proof}{\prepf\rm}{\endprepf}
\newcommand{\qed}{\hfill$\Box$}

\begin{document}

\title{Tight quasi-symmetric designs}
\author{Peter J. Cameron\\\small School of Mathematics and Statistics,
University of St Andrews\\\small North Haugh, St Andrews, Fife KY16 9SS, U.K.\\
\small email: \texttt{pjc20@st-andrews.ac.uk}}
\date{}
\maketitle

\begin{abstract}
I give an elementary proof that a quasi-symmetric design without repeated
blocks on $v$ points has at most $v\choose2$ blocks, with equality if and
only if it is a tight $4$-design.
\end{abstract}

Let $X$ be a set of $v$ points, and $\mathcal{B}$ a set of $k$-element
subsets of $X$. The \emph{degree} of $\mathcal{B}$ is the number of
cardinalities of pairwise intersections of sets in $\mathcal{B}$. The
pair $(X,\mathcal{B})$ is a \emph{$t$-design} if the number of members of
$\mathcal{B}$ containing a given $t$-subset of $X$ is a non-zero constant;
the \emph{strength} of $\mathcal{B}$ is the maximum $t$ for which
$(X,\mathcal{B})$ is a $t$-design.

In his groundbreaking thesis~\cite{delsarte}, Philippe Delsarte proved the
following result (Theorem~5.21):

\begin{theorem}
With the above notation, suppose that $\mathcal{B}$ has degree $s$ and
strength $t$; put $e=\lfloor t/2\rfloor$. Then
\[{v\choose e}\le|\mathcal{B}|\le{v\choose s}.\]
If equality holds in one bound, then it holds in the other.
\end{theorem}

His actual result is much more general, applying to arbitrary Q-polynomial
association schemes, defined for the first time in his thesis; the proof
uses the linear programming bound for subsets of association schemes, also
introduced there.

Here I am going to give a self-contained proof of a much weaker result.
A \emph{quasi-symmetric design} is a $2$-design of degree~$2$. See~\cite{ss}
for detailed information about these designs.

\begin{prop}
A quasi-symmetric design on $v$ points has at most $v\choose 2$ blocks, with
equality if and only if it is a $4$-design.
\end{prop}

Delsarte called a design meeting his bound \emph{tight}. The tight
$4$-designs were subsequently determined by Enomoto \emph{et al.}~\cite{ein}
and Bremner~\cite{bremner}; apart from the trivial example where $\mathcal{B}$
consists of all $(v-2)$-subsets of $X$, the only examples are the Witt
design (with $v=23$, $k=7$) and its complement.

\begin{proof}
Let $(X,\mathcal{B})$ be a quasi-symmetric design with $|X|=v$ and 
$|B|=k$ for $B\in\mathcal{B}$. Let $|B_1\cap B_2|\in\{x,y\}$ for all
distinct $B_1,B_2\in\mathcal{B}$, with $x<y$. The \emph{block graph} of
the design is the graph $\Gamma$ whose vertex set is $\mathcal{B}$, with $B_1$
and $B_2$ adjacent if $|B_1\cap B_2|=y$. Goethals and Seidel~\cite{gs} showed
that this graph is strongly regular. Let $A$ be its adjacency matrix.

Also, let $T(v)$ be the \emph{triangular graph} whose vertex set is the set
of $2$-subsets of $X$, two pairs adjacent if they intersect in one point.

Let $M$ be the incidence matrix of point pairs and blocks. That is, the rows of
$M$ are indexed by point pairs, and the columns by blocks; the
$(\{x,y\},B)$ entry is $1$ if $\{x,y\}\subseteq B$, and $0$ otherwise.

The strategy of the proof is as follows:
\begin{enumerate}\itemsep0pt
\item Show that $M$ has rank equal to $|\mathcal{B}|$; it follows that
$|\mathcal{B}|\le{v\choose 2}$.
\item Assuming that equality holds, show that the eigenspaces of $MM^\top$
coincide with the those of the adjacency matrix of the triangular graph. It
follows that $MM^\top$ is a linear combination of the identity matrix, the
all-$1$ matrix, and the adjacency matrix of $T(v)$. This implies that the
$(\{x,y\},\{z,u\})$ entry of $MM^\top$ is constant over all pairs with
$\{x,y\}\cap\{z,u\}=\emptyset$; in other words, any four points lie in a
constant number of blocks.
\end{enumerate}
The first step requires a short computation, while the second has a small
subtlety.

\subparagraph{Proof of (a)}
Standard arguments show that the eigenvalues of $\Gamma$ (apart from the
valency) are
$(r-\lambda-k-x)/(y-x)$ with multiplicity~$v-1$, and $-(k-x)/(y-x)$ with
multiplicity $b-v$.

Now let $M$ be as defined above. We easily find that
\[M^\top M={k\choose 2}I+{y\choose 2}A+{x\choose 2}(J-I-A).\]
From this it is easy to compute the eigenvalues of $M^\top M$ and show that
it is invertible. We don't have to keep track of the eigenvalue on the all-$1$
eigenvector, which is certainly positive; so I won't calculate the coefficient
of $J$ in the equations below. We have
\[M^\top M=(k-x)(k+x-1)/2)I+((y-x)(y+x-1)/2)A+(*)J,\]
so the eigenvalues of $M^\top M$ on vectors orthogonal to $j$ are of the
form $(k-x)(k+x-1)+(y-x)(y+x-1)\alpha$, where $\alpha$ is an eigenvalue
of $A$. This is certainly non-zero for the positive value of $\alpha$.
So put $\alpha=-(k-x)/(y-x)$, and compute the eigenvalue of $M^\top M$ to be
\[\frac{(k-x)(k+x-1)}{2}-\frac{(y-x)(y+x-1)}{2}\cdot\frac{(k-x)}{(y-x)}
=\frac{(k-x)(k-y)}{2}\ne0.\]

\subparagraph{Proof of (b)} Assume that $|\mathcal{B}|={v\choose 2}$. Then
$M$ is an isomorphism from $\mathbb{R}^{v\choose 2}$ to
$\mathbb{R}^{\mathcal{B}}$. So $MM^\top$ and $M^\top M$ have the same
spectrum, and $M^\top$ maps eigenspaces of $M^\top M$ to eigenspaces of
$MM^\top$. We know that the former coincide with those of the block graph
$\Gamma$, and need to show that the latter coincide with those of $T(v)$.

The adjacency matrix of the triangular graph has
three eigenspaces on $\mathbb{R}^{v\choose 2}$. The first is spanned by the
all-$1$ vector, and the sum of the first two is spanned by all vectors 
(regarded as functions on the set of $2$-subsets) of the form $f_x$, where
$f_x(\{y,z\})=1$ if $x\in\{y,z\}$, $0$ otherwise. Similarly, the adjacency
matrix of the block graph has three eigenspaces on $\mathbb{R}^{\mathcal{B}}$;
the first is spanned by the all-$1$ vector, and the sum of the first two by
vectors of the form $g_x$, where $g_x(B)=1$ if $x\in B$, $0$ otherwise.

Now it is clear that $M$ maps the first eigenspace of $T(v)$ to the first
eigenspace of $\Gamma$, and $M^\top$ does the reverse. Also, since
$f_xM=(k-1)g_x$, $M$ maps the sum of the first two eigenspaces of $T(v)$
to the sum of the first two eigenspaces of
$\Gamma$. Hence, by duality, $M^\top$ maps the third eigenspace of $\Gamma$
(the orthogonal complement of the sum of the first two) to the third
eigenspace of $T(v)$; since $MM^\top$ is symmetric, the sum of the other two
eigenspaces of $MM^\top$ is the orthogonal complement of this space, hence
the sum of the first two eigenspaces of $T(v)$. \qed
\end{proof}

\paragraph{Acknowledgement} I am grateful to Bhaskar Bagchi for encouragement
and helpful comments on this note.

\end{document}